\documentclass[11pt,a4paper]{amsart}

\usepackage{newtxtext}       %

\usepackage{amsmath}
\usepackage{amssymb}
\usepackage{MnSymbol}
\usepackage{bm}
\usepackage{hyperref}
\usepackage{url}

\usepackage{cleveref}

\DeclareMathOperator{\Expectation}{\mathbb E}
\DeclareMathOperator{\Maxexp}{\mathcal E}

\newcommand{\KH}[2]{\operatorname{DH}\left(#1 \middle| #2 \right)}

\newcommand{\avalof}[1]{\left\vert#1\right\vert}

\newcommand{\coshtwo}{\cosh_2}

\newcommand{\cpoly}[2]{C^{#1}_{\text{poly}}\left(\reals^{#2}\right)}
\newcommand{\derivby}[1]{\frac{d}{d#1}}
\newcommand{\euler}{\mathrm{e}}
\newcommand{\expectat}[2]{\Expectation_{#1}\left[#2\right]}

\newcommand{\expof}[1]{\exp\left(#1\right)}

\newcommand{\gaussdensity}{\gamma}
\newcommand{\gaussint}[2]{\int{#1}\ \gaussdensity(#2) \ d#2 \ }

\newcommand{\maxexpat}[1]{\Maxexp\left(#1\right)}
\newcommand{\naturals}{\mathbb{N}}
\newcommand{\normat}[2]{\left\Vert#2\right\Vert_{#1}}

\newcommand{\orliczof}[2]{L_{#1}\left(#2\right)}
\newcommand{\orliczpof}[3]{L_{#1}^{#2}\left(#3\right)}
\newcommand{\pderivby}[1]{\frac {\partial} {\partial #1}}

\newcommand{\reals}{\mathbb{R}}
\newcommand{\scalarat}[3]{\left\langle#2,#3\right\rangle_{#1}}

\newcommand{\setof}[2]{\left\{#1 \, \middle| \, #2 \right\}}
\newcommand{\set}[1]{\left\{#1\right\}}

\begin{document}

\title[Orlicz Spaces in IG]{A lecture about the use of Orlicz Spaces
  in Information Geometry}
\author[G. Pistone]{Giovanni Pistone}
\address{de Castro Statistics, Collegio Carlo
  Alberto, Torino, Italy}
\email{giovanni.pistone@carloalberto.org}
\urladdr{https;//www.giannidiorestino.it}

\maketitle

\begin{abstract}This chapter is a revised version of a tutorial lecture that
  I presented at the \'Ecole de Physique des Houches on July 26-31
  2020. Topics include: Non-parametric Information Geometry, the
  Statistical bundle, exponential Orlicz spaces, and Gaussian
  Orlicz-Sobolev spaces
\end{abstract}

\section{Introduction}
\label{sec:introduction}
This chapter is a revision of the lecture and the related hand-out
which I presented to the École de Physique des Houches on July 26-31
2020. Due to its strictly tutorial character, I shall not give
detailed primary references. I shall mention some references that
expand and support specific points, and I shall add some in a final
section to point to further developments.

I aim to review the basics of a peculiar setting for Information
Geometry (IG) in the sense of \cite{amari|nagaoka:2000}. This setting
has the following peculiarities.

\begin{itemize} \item It is non-parametric and infinite-dimensional.
\item It provides an affine manifold modeled on a Banach space in the sense of \cite{lang:1995}, the Banach space being an Orlicz space as defined \cite[Ch. 8]{adams|fournier:2003}.
\item  It focuses on a particular expression of the tangent bundle, called Statistical Bundle (ST).
\item It allows for the use of (weakly) differentiable densities when the reference measure is Gaussian as in \cite[Ch. V]{malliavin:1995}.
\end{itemize}

A previous tutorial paper \cite{pistone:2020NPCS} presents this
non-parametric construction in the case of finite state space. A
general presentation of the finite case is in
\cite[Ch. 2]{Ay|Jost|Le|Schwachhofer:2017IGbook}. Here, I will focus
on the preliminaries of the infinite state space case. A comprehensive
modern introduction to the whole topic is found in
\cite{nielsen|2020-entropy-1100}.

There are many other successful presentations of IG which are indeed
non-parametric. I will give a few relevant references in the concluding
section. Any useful presentation should explain and
include its historical development elements
below $\bm 1$ to $\bm 5$.

\noindent $\bm 1.$ The starting point to consider is the now classical
work of R. Fisher, see, for example,
\cite[Ch. 4]{efron|hastie:2016}. A \emph{regular statistical model} is
a mapping from a set of parameters $\Theta$ which is an open domain of
$\reals^d$ to probability densities on a given measured sample space
$\mathcal P(X,\mathcal X,\mu)$, $\theta \mapsto p(\theta)$, such that
the following computation is feasible. If $f$ is a given random
variable, one wants to compute the variation of the expectation
(assuming its existence), as
\begin{multline*}
  \pderivby {\theta_i} \expectat {p(\theta)} f = \pderivby {\theta_i}
  \int f\ p(\theta) \ d\mu = \pderivby {\theta_i} \scalarat \mu f
  {p(\theta)} = \\
  \scalarat \mu f {\pderivby {\theta_i} p(\theta)} = \scalarat \mu f
  {\frac{\pderivby {\theta_i} p(\theta)}{p(\theta)} p(\theta)} = \scalarat {p(\theta)} f
  {\pderivby {\theta_i} \log p(\theta)} = \\
  \scalarat {p(\theta)} {f - \expectat {p(\theta)} f } 
  {\pderivby {\theta_i} \log p(\theta)} \ .
\end{multline*}

The random vector with components $\pderivby {\theta_i} \log
p(\theta)$ is the \emph{Fisher score} of the model al $\theta$. Its
expected value with respect to $p(\theta)$ is 0 and its
variance matrix with respect to $p(\theta)$ is the \emph{Fisher
information matrix},
\begin{equation}\label{eq:fisher}
  I(\theta) = \left[ \expectat {p(\theta)} {\pderivby {\theta_i} \log
      p(\theta) \pderivby {\theta_j} \log p(\theta)} \right]_{ij} =
  \left[ \int \frac{\pderivby {\theta_i}  p(\theta) \pderivby
      {\theta_j}  p(\theta)}{p(\theta)} \ d\mu \right]_{ij} \ .
\end{equation}

The computations above shows many peculiar  features:
\begin{itemize} \item  One computes the variation of the expectation
  as a function of the statistical model. Moreover, a moving inner
  product $\theta \mapsto \scalarat {p(\theta)} \cdot \cdot$ appears naturally.
\item The score $\partial_i \log p(\theta)$ represents the velocity of
  variation of the statistical model, while $f - \expectat {p(\theta)}
  f$ represents the gradient of the expectation function.
\item  The velocity at $\theta$ lives in the space of random variables
  centred at $p(\theta)$.
\item The information matrix provides squared norms and scalar
  products of the velocities in the moving inner product.
\end{itemize}

\noindent$\bm 2.$ One explanation of the Fisher computations results
from the assumption of an exponential model,
\begin{equation*}
  p(\theta) = \euler^{\sum_i \theta_i u_i - \kappa(\theta)} \ ,
\end{equation*}
where $u_i$ are the sufficient statistics and $\kappa(\theta)$ is the
cumulant of $\sum_i \theta_i u_i$, see, for example,
\cite[Ch. 5]{efron|hastie:2016}. In such a case, the score is the
centered sufficient statistics,
\begin{equation*}
  \partial _i \log p(\theta) = u_i - \partial_i
  \kappa(\theta) = u_i - \expectat {p(\theta)} {u_i} \ ,
\end{equation*}
and the information matrix equals the Hessian of the convex function
$\theta \mapsto \kappa(\theta)$,
\begin{equation*}
  I(\theta) = \left[\partial_i \partial_j \kappa(\theta)\right]_{ij} \ .
\end{equation*}
See \cite[Ch. 2]{brown:86} for a full account of analytic properties
of exponential families.

\noindent$\bm 3.$ C. R. Rao has been the first statistician to remark
that the Fisher information matrix is positive definite and smooth in
an adequately defined regular model. Hence it defines a Riemannian
metric on the space of parameters. Moreover, he provided an embedding
argument for the resulting manifold. What possibly inspired him in his
construction was Differential Geometry use in Physics. The original
presentation of IG in \cite{amari|nagaoka:2000} can be considered the
full unfolding of this approach.

The mapping
\begin{equation*}
  \Theta \ni \theta \mapsto 2 \sqrt{p(\theta)} = P(\theta)
\end{equation*}
maps the parameters' space $\Theta$ into the $L^2(\mu)$-sphere of
radius $2$. The vectors
\begin{equation*}
  \partial_i P(\theta) = \partial_i 2 \sqrt{p(\theta)} =
  \frac{\partial_i p(\theta)}{\sqrt{p(\theta)}}
\end{equation*}
are in the tangent space at $P(\theta)$ of the sphere, and the inner
product between tangent vector is
\begin{equation*}
  \int \partial_i P(\theta) \partial_j P(\theta) \ d\mu = \int
  \frac{\partial_i p(\theta) \partial_j(x,\theta)}{p(x;\theta)} \ d\mu \ ,
\end{equation*}
that is, the $(i,j)$ element of the Fisher information matrix.

The Rao's computations above reproduce all the metric structure of the
Fisher computations but in one point. That is, now the velocity is not
expressed by the logarithmic derivative $\partial_i \log p(\theta) =
\partial_i p(\theta) / p(\theta)$, but it is expressed by $\partial_i
2 \sqrt{p(\theta)} = \partial_i p(\theta) / \sqrt{p(\theta)}$.

It is possible to make sense, at least formally, of the apparent
contradiction by considering that there are here three different
expressions of the same object.
\begin{itemize}
\item $T\mathcal P_>$ is the tangent bundle of the set of positive
  densities $\mathcal P_>$; the tangent vectors $\dot p$ satisfy
  $\int \dot p \ d\mu = 0$.
\item $TS_2$ is the tangent bundle of the $L^2(\mu)$ sphere $S_2$.
\item $S \mathcal P_>$ is the Fisher's statistical bundle consisting
  of all couples $(p.u)$ such that $p \in \mathcal P_>$ and
  $\expectat p u = 0$.
\end{itemize}

The Rao's embedding $p \mapsto 2 \sqrt p$ provides the identification
of $T\mathcal P$ with $TS_2$. In fact, the computation of the tangent
mapping of $s \colon P \mapsto P^2/4 = p$ gives $ds(P)[\dot P] = \sqrt
p \dot P$.

The identification of $TS_2$ with
$S\mathcal P$ is provided by
\begin{equation*}
  TS_2 \ni (P,\dot P) \mapsto \left(\frac14 P^2, 2 \frac {\dot P}
    P\right) = (p,u) \in S\mathcal P \ .
\end{equation*}
In fact,
\begin{gather*}
  \int p \ d\mu = \frac14 \int P^2 \ d\mu = 1 \ ,\\
  \int u \ p \ d\mu = \frac 12 \int \frac {\dot P}{P} P^2 \ d\mu = \frac
  12 \int P \dot P \ d\mu = 0 \ ,\\
  \int u_1 u_2 \ p \ d\mu = \int \frac{\dot P_1}{P} \frac{\dot P_2}{P}
  P^2 \ d\mu = \int \dot P_1 \dot P_2 \ d\mu \ . 
\end{gather*}

A large part of the literature in IG uses the expressions
$T\mathcal P$ and $TS_2$. Still, my own choice is to use
$S\mathcal P$.  It fits well with the statistical picture and the
exponential representation of strictly positive densities.

The choice of the exponential expression and Fisher's score might seem
arbitrary, but the following argument shows is not, see
\cite[Ch. 3]{Ay|Jost|Le|Schwachhofer:2017IGbook}. Assume
$t \mapsto \mu(t)$ is a one-dimensional smooth model of probability
measures and assume the mapping $t \mapsto \mu(A;t)$ is smooth for all
measurable $A$. Fix a value $\bar t$ of the parameter $t$. If a
measurable set $A$ is a zero set for $\mu(\bar t)$, then
$t \mapsto \mu(A,t)$ is minimum at $t = \bar t$ then the derivative is
zero at $\bar t$, namely $\dot \mu(\bar t)=0$. Il follows that the
measure $\dot \mu(\bar t)$ is absolutely continuous with respect to
$\mu(\bar t)$. The resulting logarithmic derivative
$d \dot \mu(\bar t) / d \mu(\bar t)$ is clearly a generalization of
the Fisher's score.

\noindent$\bm 4.$ The Riemannian approach by C. R. Rao can lead to a
more in-depth study of the second-order properties of the manifold,
namely the Levi-Civita connection and the curvature. Several authors,
notably S-I. Amari, B. Efron, Ph. Dawid, S. Lauritzen have promoted a more general perspective in studying the geometry of statistical
models. In modern terminology, a
statistical manifold consists of a metric and a couple of flat
connections,  in duality for the given inner product. This
set-up nicely solves the divide between Fisher's and Rao's
approach by producing a unified theory. Moreover, the specific type of
affine manifold relevant for IG is a Hessian manifold. That
is, all its structure depends on a master convex functional.

\noindent$\bm 5.$ Both theoretical and applied research have recently
shown interest in a particular type of non-parametric statistical
models. Namely, models where the real space $\reals^n$ is a model for
the sample space, the reference measure is either the Lebesgue measure
or the Gaussian measure and the densities are required to have some
level of smoothness.

For example, there is the statistical estimation method based on
Hyv\"arinen's divergence,
\begin{equation} \label{eq:hyvarinen}
  \KH p q = \frac12 \int \avalof{\nabla \log p(x) - \nabla \log q(x)}^2 \ p(x) \ dx \ , 
\end{equation}
see \cite[\S 13.6.2]{amari:2016}, 
where $p,q$ are positive probability densities of the $n$-dimensional
Lebesgue space. It is assumed that the $\log$-densities are smooth and
the integral exists. If we define the Otto's inner product by
\begin{equation} \label{eq:otto}
  \llangle f,g \rrangle _ p = \int \nabla f(x) \cdot \nabla g(x) \ p(x) \ dx \ ,
\end{equation}
where $p$ is a probability density and $f,g$ are smooth random
variables such that $\expectat p f = \expectat p g = 0$, then the
development of the square in the Hyv\"arinen divergence produces the
term
\begin{equation*}
  \llangle \log p - \expectat p {\log p} , \log q - \expectat p {\log q}
  \rrangle _ p \ . 
\end{equation*}
Notice that the equation above with $p=q$ reminds of the Fisher
information of \cref{eq:fisher}, with parameter derivative replaced by
spatial derivatives. In particular, parameter and spatial derivatives
coincide in the case of translation models, that is, when
$p(x,\theta) = p(x+\theta)$.

In the Gaussian case, the celebrated log-Sobolev inequality, see
\cite[Ch. 5]{ledoux:2001-concentration}, can be written as
\begin{multline*}
  \operatorname{Ent}(p) = \int p(x) \log p(x) \ \gamma(x) dx \leq \\ 2
  \int \avalof {\nabla \sqrt {p(x)}} ^2 \ \gamma(x)dx = \frac 12
\llangle \log p - \operatorname{Ent}(p) , \log p - \operatorname{Ent}(p) 
  \rrangle _ p \ , \end{multline*}
where $\gamma$ is the standard Gaussian density, $p$ is a density with
respect to $\gamma$ and $\avalof \cdot$ is the Euclidean norm.

In the following, I shall focus on the exponential representation of
positive densities $p = \euler^{u - K(u)}$. The reference measure is
specialised to be $\mu(dx) = \gaussdensity(x) dx$ where
$\gaussdensity$ is the standard Gaussian density of $\reals^n$. The
sufficient statistics $u$ is assumed to belong to an exponential
Orlicz space, to be defined next, and such that
$\gaussint {u(x)} x = 0$. The normalizing constant $u \mapsto K(u)$ is
a convex function defined on the exponential Orlicz space. The interior of
the proper domain of $K$ is the parameter set of a maximal
exponential Gaussian space model. Here, maximal means that the
model contains all possible finite dimensional exponential families. I
will aim to sketch a theory in which all concerns of items 1 to 5 meet
a solution of a sort.

\section{Orlicz spaces}
\label{sec:stat-grad-model}
First, let us review a few elements of the theory of Orlicz spaces and
fix a convenient notation. I do not not aim to full generality,
cf. \cite[Ch. 8]{adams|fournier:2003}. If $\phi \in C[0,+\infty[$ satisfies:
\begin{enumerate}
\item  $\phi(0)=0$,
\item $\phi$ is strictly increasing, and
\item $\lim_{u \to +\infty} \phi(u) = +\infty$,
\end{enumerate}
then its primitive function
\begin{equation*}
  \Phi(x) = \int_0^x \phi(u) \ du \ , \quad x \geq 0 \ ,
\end{equation*}
is strictly convex. The function $\Phi$ is extended to $\reals$ by
symmetry, $\Phi(x) = \Phi(\avalof x)$, and is called \emph{Young
  function}.

The inverse function $\psi = \phi^{-1}$ has the same properties 1) to 3) as $\phi$, so that its primitive
\begin{equation*}
  \Psi(y) = \int_0^y \psi(v) \ dv \ , \quad y \geq 0 \ ,
\end{equation*}
is again a Young function. The couple $(\Phi, \Psi)$, is a couple of \emph{conjugate} Young functions. The relation is symmetric and we write both $\Psi=\Phi_*$ and $\Phi = \Psi_*$. The Young inequality holds true,
\begin{equation*}
  \Phi(x) + \Psi(y) \geq xy \ , \quad x,y \geq 0 \ ,
\end{equation*}
and the Legendre equality holds true ,
\begin{equation*}
  \Phi(x) + \Psi(\phi(x)) = x \phi(x) \ , \quad x \geq 0 \ .
\end{equation*}

Here are the specific cases I am going to use. The sub-2 index
denotes the 2nd Taylor remainder.
\begin{align}
&\Phi_\alpha(x) = \frac {x^\alpha} \alpha \ , \quad \Psi_\beta(y) =
                \frac {y^\beta}\beta \ , \quad \alpha, \beta > 1 \  , \quad \frac1\alpha+\frac1\beta = 1 \ ;   \label{eq:cases-1} \\
  &\exp_2(x) = \euler^x - 1 - x \ , \quad (\exp_2)_*(y) = (1+y)\log(1+y) - y \ ; \label{eq:cases-2} \\
  &\cosh_2(x) = \cosh x - 1 \ , \quad (\cosh_2)_*(y) = \int_0^y \sinh^{-1}(v) \ dv \ ; \label{eq:cases-3} \\
&\operatorname{gauss}_2(x) = \expof{\frac12x^2}-1 \ . \label{eq:cases-4}
\end{align}

Given a Young function $\Phi$ and a probability measure $\mu$, the
\emph{Orlicz space} $\orliczof \Phi \mu$ is the Banach space whose closed
unit ball is $\setof{f \in L^0(\mu)}{\int \Phi(\avalof f) \ d\mu \leq
  1}$. This defines the \emph{Luxemburg norm}, characterized by
\begin{equation*}
  \normat {\orliczof \Phi \mu} f \leq \rho \quad \text{if, and only if,} \quad \int \Phi(\rho^{-1} \avalof f) \ d\mu   \leq 1 \ .
\end{equation*}

Because of the Young inequality, it holds
\begin{equation*}
  \int \avalof{uv} \ d\mu \leq \int \Phi(\avalof u) \ d\mu + \int \Phi_*(\avalof v) \ d\mu \ .
\end{equation*}
This provides a separating duality $\scalarat \mu u v = \int uv \ d\mu$ of $\orliczof {\Phi} \mu$ and $\orliczof {\Phi_*} \mu$ such that
\begin{equation*}
  \scalarat \mu u v \leq 2 \normat {\orliczof{\Phi} \mu} u  \normat {\orliczof {\Phi_*} \mu} v \ .
\end{equation*}
From the conjugation between $\Phi$ and $\Psi$, an equivalent norm can
be defined, namely, the \emph{Orlicz norm}
\begin{equation*}
  \normat {{\orliczof {\Phi} \mu }^*} f = \sup \setof{\scalarat \mu f
    g}{\normat {\orliczof {\Psi} \mu} f \leq 1} \ .
\end{equation*}

The domination relation between Young functions imply continuous
injection properties for the corresponding Orlicz spaces. We wll say
that $\Phi_2$ \emph{eventually dominates} $\Phi_1$, written
$\Phi_1 \prec \Phi_2$, if there is a constant $k$ such that
$\Phi_1(x) \leq \Phi_2(kx)$ for all $x$ larger than some $\bar x$. As,
in our case, $\mu$ is a probability measure, the continuous embedding
$\orliczof {\Phi_2} \mu \to \orliczof {\Phi_1} \mu$ holds if, and only
if, $\Phi_1 \prec \Phi_2$. If $\Phi_1 \prec \Phi_2$, then
$(\Phi_2)_* \prec (\Phi_1)_*$.

A special case occurs when there exists a function $C$ such that
$\Phi(ax) \leq C(a) \Phi(x)$ for all $a \geq 0$, which is the case, for
example, for a power function and in the case of the functions
$(\exp_2)_*$ and $(\cosh-1)_*$. In such a case, the dual couple is a
couple of reflexive Banach spaces, and bounded functions are a dense
set. I will return to this important topic below.

The spaces corresponding to power case
\eqref{eq:cases-1} coincides with the  ordinary
Lebesgue spaces. The norm are related by
\begin{equation*}
  \normat {\orliczof {\Phi_\alpha} \mu } f = \alpha^{1/\alpha} \normat
  {L^\alpha(\mu)} f \ .
\end{equation*}
With reference to our examples \eqref{eq:cases-2} and
\eqref{eq:cases-3}, we see that $\exp_2$ and $(\cosh-1)$ are
equivalent. They both are eventually dominated by
$\operatorname{gauss}_2$ \eqref{eq:cases-4} and eventually dominate
all powers \eqref{eq:cases-1}.  The cases \eqref{eq:cases-2} and
\eqref{eq:cases-3} provide isomorphic B-spaces
$\orliczof {(\cosh-1)} \mu \leftrightarrow \orliczof {\exp_2} \mu$
which are of special interest for us as they provide the model spaces
for our non-parametric version of IG, see
\cref{sec:exponential-bundle} below.

Clearly, a function belongs to the space $\orliczof {\cosh_2} \mu$ if,
and only if, its \emph{moment generating function
$\lambda \mapsto \int \euler^{\lambda f}$ is finite in a neighborhood
of 0}. In turn, this implies that the moment generating function is
analytic at 0, see \cite[Ch. 2]{brown:86}.

The same property is equivalent to a large deviation inequality, see
\cite[Ch. 2]{wainwright:2019-HDS}. \emph{A function $f$ belongs to
  $\orliczof {\cosh-1} \mu$ if, and only if, it is
  \emph{sub-exponential}, that is, there exist constants $C_1,C_2 > 0$
  such that}
\begin{equation*}
  \mu (\avalof{f} \geq t) \leq C_1 \expof{-C_2 t} \ , \quad t \geq 0 \ .
\end{equation*}

Let us check the equivalence above. If $\normat {\orliczof {\cosh_2}
  \mu} f = \rho$, then $\int \euler^{\rho^{-1} \avalof f} \ d\mu \leq 4$. It
  follows that
  \begin{equation*}
    \mu(\avalof f > t) =
    \mu\left(\euler^{\rho^{-1}\avalof f} >
      \euler^{\rho^{-1}t}\right) \leq
    \left(\int \euler^{\rho^{-1}\avalof f} \ d\mu\right) 
    \euler^{- \rho^{-1}t} \leq
    4 \euler^{- \rho^{-1}t} \ .  
  \end{equation*}
The sub-exponential inequality holds with $C_1=4$ and $C_2 = \normat
{\orliczof {\cosh_2} \mu} f ^{-1}$. Conversely, for all $\lambda > 0$,
\begin{equation*}
  \int \euler^{\lambda f} \ d\mu \leq \int_1^\infty \mu\left(\euler^{\lambda
      f^+}> t\right) \ dt \leq C_1 \int_0^\infty \euler^{-(C_2
    \lambda^{-1} - 1)s} \ ds \ .  
\end{equation*}
The right-hand side is finite if $\lambda < C_2$ and the same
bound holds for $-f$.

A sub-exponential random variable is of particular interest in
applications because they admit an explicit exponential bound in the
Law of Large Numbers. Another class of interest consists of
the \emph{sub-Gaussian} random variables, that is, those
random variables whose square is sub-exponential.

The theory of sub-exponential random variables provides an
\emph{equivalent norm for the space} $\orliczof {\cosh_2} \mu$, namely
the norm
\begin{equation*}
  f \mapsto \sup_k \left((2k)!^{-1} \int f^{2k} \ d\mu
  \right)^{1/2k} = \left\bracevert f \right\bracevert _ {\cosh_2} \ .
\end{equation*}
See \cite{buldygin|kozachenko:2000} or
\cite{siri|trivellato:2020-SPL}. Let us prove the equivalence. If
$\normat {\orliczof {\cosh_2} \mu} f \leq 1$, then
\begin{equation*}
  1 \geq \int \cosh_2 f \ d\mu \geq \frac1{(2k)!} \int f^{2k} \ d\mu
    \quad \text{for all $k = 1,2,\dots$,}
\end{equation*}
so that $1 \geq  \left\bracevert f \right\bracevert _
{\cosh_2}$. Conversely, if the latter inequality holds, then
\begin{equation*}
  \int \coshtwo (f/\sqrt 2) \ d\mu = \sum_{k=1}^\infty \frac1{(2k)!}
  \int f^{2k} \ d\mu \left(\frac12\right)^k \leq 1 \ , 
\end{equation*}
so that $\normat {\orliczof {\coshtwo} \mu} f \leq \sqrt 2$.

It is covenient to introduce a further notation. For each Young
function $\Phi$, the function $\overline \Phi(x) = \Phi(x^2)$ is again
a Young function such that
$\normat {\orliczof {\overline \Phi} \mu} f \leq \lambda$ if, and only
if, $\normat {\orliczof \Phi \mu} {\avalof f ^2} \leq
\lambda^2$. \emph{We will denote the resulting space by}
$\orliczpof \Phi 2 \mu$. For example, $\operatorname{gauss}_2$ and
$\overline{\cosh-1}$ are $\prec$-equivalent , hence the isomorphisn
$\orliczof {\operatorname{gauss}_2} \mu \leftrightarrow \orliczpof
{(\cosh-1)} 2 \mu$. As an application of this notation, consider that
for each increasing convex $\Phi$ it holds
$\Phi(fg) \leq \Phi((f^2+g^2)/2) \leq (\Phi(f^2) + \Phi(g^2))/2$. It
follows that when the $\orliczpof \Phi 2 \mu$-norm of $f$ and of $g$
is bounded by one, the $\orliczof \Phi \mu$-norm of $f$, $g$, and
$fg$, are all bounded by one. The space $\orliczof {\coshtwo} \mu$ has
a continuous injection in the Fr\'echet space
$L^{\infty-0}(\mu) = \cap_{\alpha>1} L^\alpha(\mu)$, which is an
algebra. When we need the product, we can either assume the factor are
both sub-Gaussian, or, move up the functional frawork to the
intersection of the Lebesgue spaces.

Let us now discuss specific issues of the Gaussian exponential Orlicz
spaces $\orliczof {\cosh_2} \gaussdensity$, $\gaussdensity$ the
standard $n$-variate Gaussian density. \emph{Dominated convergence
  does not hold in this space}. The squared-norm function
$f(x) = \avalof x ^ 2$ belongs to the Gaussian exponential Orlicz
space $\orliczof {\cosh_2} \gaussdensity$ because
\begin{equation*}
  \int \cosh_2(\lambda f(x)) \ \gaussdensity(x)dx < \infty \quad
  \text{for all} \quad \lambda < 1/2 \ .
\end{equation*}
The sequence $f_N(x) = f(x)(\avalof x \leq N)$ converges to $f$
point-wise and in all $L^\alpha(\gamma)$, $1 \leq \alpha < \infty$.
However, the convergence does not hold in the Gaussian exponential
Orlicz space. In fact, for all $\lambda \geq 1/2$,
\begin{multline*}
 \int  \cosh_2(\lambda (f(x)-f_N(x))) \ \gaussdensity(x)dx =
 \int_{\avalof x > N} \coshtwo (\lambda f(x)) \ \gaussdensity(x)dx =
 \infty \ ,
\end{multline*}
but convergence would imply
\begin{equation*}
  \limsup_{N \to \infty} \int  \cosh_2(\lambda (f(x)-f_N(x))) \ \gaussdensity(x)dx \leq 1 \quad \text{for all $\lambda >0$}\ .
\end{equation*}

The closure in $\orliczof {\cosh_2} \gaussdensity$ of the vector space
of bounded functions is called \emph{Orlicz class} and it is denoted
by $M_{\cosh_2}(\gaussdensity)$. One can prove that
$f \in M_{\cosh_2}(\gaussdensity)$ if, and only if. the moment
generating function
$\lambda \mapsto \int \euler^{\lambda f(x)} \ \gaussdensity(x)dx$ is
finite for all $\lambda$, see \cite{pistone:2018}. An example is
$f(x) = x$. Bounded confergence holds in the Orlicz class. Assume
$f \in M_{\cosh_2}(\gaussdensity)$ and consider the sequence
$f_N(x) = (\avalof x \leq N) f(x)$. Now,
\begin{equation*}
  \int \coshtwo(\lambda(f(x) - f_N(x)) \ \gaussdensity(x) dx =
  \int_{\avalof x \geq N} \coshtwo(\lambda f(x)) \ \gaussdensity(x) dz
\end{equation*}
converges to 0 as $N \to \infty$.

\section{Calculus of the Gaussian space}
\label{sec:calc-gauss-space}
I will review here a few simple facts about the analysis of the
Gaussian space, the so-called Malliavin's calculus, see
\cite[Ch.~V]{malliavin:1995}. 

Let us denote by $C^k_\text{poly}(\reals^n)$, $k = 0,1,\dots$, the vector space of functions which are differentiable up to order $k$ and which are bounded, together with all derivatives, by a polynomial. This class of functions is dense in $L^2(\gaussdensity)$. For each couple $f,g \in \cpoly 1 n$, we have
\begin{equation*}
  \gaussint {f(x) \ \partial_i g(x)} x =   \gaussint {\delta_i f(x) \ g(x)} x \ ,
\end{equation*}
where the divergence operator $\delta_i$ is defined by $\delta_i f(x) = x_i f(x) - \partial_i f(x)$. Multidimensional notations will be used, for example,
\begin{equation*}
  \gaussint {\nabla f(x) \cdot \nabla g(x)} x = \gaussint {f(x) \ \delta \cdot \nabla g(x)} x \ , \quad f,g \in \cpoly 2 n \ ,
\end{equation*}
with $\delta \cdot \nabla g(x) = x \cdot \nabla g(x) - \Delta g(x)$.

For example, in this notation, the Hyv\"arinen divergence of \cref{eq:hyvarinen} with $P = p\cdot \gaussdensity$, $Q = q \cdot \gaussdensity$, and $p,q \in \cpoly 2 n$, becomes
\begin{equation*} 
  \frac12 \int \avalof {\nabla \log p(x) - \nabla \log q(x)} ^ 2 p(x)
  \gaussdensity(x)dx \ ,
\end{equation*}
while the Otto's inner product of \cref{eq:otto} becomes, with $P = p
\cdot \gaussdensity$ and $f,g,p  \in \cpoly 2 n$, gives
\begin{equation*}
  \gaussint {\nabla f(x) \cdot \nabla g(x) \ p(x)} x  = \gaussint { f(x) \delta \cdot \nabla (g(x) p(x))} x \ .
\end{equation*}

Hermite polynomials $H_\alpha = \delta^\alpha 1$ provide an orthogonal
basis for $L^2(\gaussdensity)$ such that
$\partial_i H_\alpha = \alpha_i H_{\alpha - e_i}$. Fourier expansion
in Hermite polynomials provides proof of the closure of both operator
$\partial_i$ and $\delta_i$ on a domain which is an Hilbert subspace
of $L^2(\gaussdensity)$. Moreover, the closure of $\partial_i$ is the
translation operator's infinitesimal generator.

\section{Exponential statistical bundle}
\label{sec:exponential-bundle}
In this section, I will very briefly review and slightly generalize my
construction of the statistical manifold as a Banach manifold modeled
on the exponential Orlicz space $\orliczof {\coshtwo}
\gaussdensity$. For a detailed presentation, see
\cite{pistone:2013GSI} and \cite{pistone:2018}. The general set-up is
specialized to the Gaussian space.

The support of the manifold is the \emph{maximal exponential model}
$\maxexpat \gaussdensity$ consisting of probability densities on
$(\reals^n,\gaussdensity)$ of the form
\begin{equation*}
  q = \expof{u - K_1(u)} \ , \quad u \in B_1 = \orliczof {\coshtwo}
  \gaussdensity \cap \setof{u}{\int u(x) \ \gaussdensity(x)dx = 0} \ .
\end{equation*}
The quantity $K_1(u) = \log \int \euler^{u(x)} \ \gaussdensity(x) dx$
is the unique normalising constant ($\log$- partition function) of
$\euler^u$ and is assumed to be finite. A further restriction is
needed to avoid the border of the set $\set{K_1(u)<\infty}$. We can
easily prove that the
mapping $K_1 \colon \orliczof {\coshtwo} \gaussdensity \to \reals$ is
convex and that the topological interior of its proper domain in
non-empty because it contains the
open unit ball of $B_1$. We restrict the model to all $u$ in such
domain. The mapping
\begin{equation*}
  s_1 \colon \maxexpat \gaussdensity \ni q \mapsto u = \log q - \expectat \gaussdensity {\log q}
  \in B_1
\end{equation*}
provides a global chart to the manifold. If we chose to express the
velocity of a curve $t \mapsto \maxexpat \gaussdensity$ by the
Fishers score $\derivby t \log q(t)$, then the tangent bundle of the
manifold is expressed by the \emph{statistical bundle}
$S \maxexpat \gaussdensity$ consisting of all the couples $(q,v)$ such
that $q$ is a density of the maximal exponential model and $v$ is a
$q$-centered random variable in the exponential Orlicz space.

The following statement is crucial to prove consistency in infinite
dimensions. It shows that the statistical bundle fibers are isomorphic
as Banach spaces.

\emph{For all $p, q \in \maxexpat \gaussdensity$ it holds $q = \euler^{u - K_p(u)} \cdot p$, where $u \in \orliczof {(\cosh-1)} \gaussdensity$, $\expectat p u = 0$, and $u$ belongs to the interior of the proper domain of the convex function $K_p$. This property is equivalent to any of the following:
\begin{enumerate}\label{portmanteaux}
\item\label{portmanteaux-1} $p$ and $q$ are connected by an open exponential arc;
\item\label{portmanteaux-2} $\orliczof {\coshtwo}
  p = \orliczof {\coshtwo}
  q$ and the norms are equivalent;
\item\label{portmanteaux-3}
  $p/q \in \cup_{a > 1} L^a(q)$ and $q/p \in \cup_{a>1} L^a(p)$.
\end{enumerate}}
See
\cite{santacroce|siri|trivellato:2016,santacroce|siri|trivellato:2018}
for a detailed proof. The following argument is a generalization of the the proof of
\cref{portmanteaux-1}$\Rightarrow$\cref{portmanteaux-2}. Let $F$ be
logarithmically convex on $\reals$ and such that $\Phi = F - 1$ is a
Young function. For example, the assumption holds for both
$F(x)=\cosh x$ and $F(x) = \euler^{x^2/2}$. For all real $A$ and $B$,
the function
\begin{equation*}
  \reals^2 \ni  (\lambda,t) \mapsto F(\lambda A) \euler^{tB} = \expof{\log F(\lambda
    A) + tB}  
\end{equation*}
is convex and so is
\begin{equation*}
  C(\lambda,t) = \int F(\lambda f(x)) \euler^{tu(x)} \
  \gaussdensity(x)dx \ ,
\end{equation*}
where $f \in \orliczof \Phi \gaussdensity$ with and
$u \in \orliczof {\cosh_2} \gaussdensity$ with
$\int u(x) \ \gaussdensity(x)dx = 0$. Without restriction of
generality, assume $\normat {\orliczof \Phi \gaussdensity} f
=1$. Let us derive two marginal inequalities. First, for $t=0$, the
definition of Luxemburg norm gives
\begin{equation*}
  C(\lambda,0) = \int F(\lambda f) \ \gaussdensity(x) dx
  \leq 2 \ , \quad -1 \leq \lambda \leq 1 \ .
\end{equation*}
Second, for $\lambda = 0$, consider
$K_1(tu) = \log \int \euler^{tu} \ \gaussdensity(x)dx$,
where $t$ belongs to an an open interval $I$ containing $[0,1]$ and
such that $K_1(tu) < + \infty$. It follows that
\begin{equation*}
  C(0,t) = \int \euler^{tu} \ \gaussdensity(x)dx = \euler^{K_1(tu)} <
  + \infty \ .
\end{equation*}

Choose a $t>1$ in $I$ and consider the convex combination
\begin{equation*}
  \left(\frac{t-1}t,1\right) = \frac{t-1}t (1,0)+\frac1t (0,t)
\end{equation*}
and the inequality
\begin{equation*}
  C\left(\frac{t-1}t,1\right) \leq \frac{t-1}t C(1,0)+\frac1t C(0,t)
  \leq 2 \frac{t-1}t + \frac1t \euler^{K_1(tu)} \ .
\end{equation*}
Now,
\begin{multline*}
  \int \Phi\left(\frac{t-1}t f(x)\right) \euler^{u(x)-K_1(u)}
  \gaussdensity(x)dx = \\ \int F\left(\frac{t-1}t f(x)\right)
  \euler^{u(x)-K_1(u)} \gaussdensity(x)dx - 1 = \\
  \euler^{-K_1(u)} C\left(\frac{t-1}t,1\right) - 1 \leq
   \euler^{-K_1(u)} \left(2\frac{t-1}t + \frac1t
     \euler^{K_1(tu)}\right) - 1 \end{multline*}

 As the right-hand-side is finite, we have proved that
 $f \in \orliczof \Phi p$ for $p = \euler^{u - K_1(u)}$. Conversely, a
 similar argument shows the other implication. We have proved that all
 Orlicz spaces $\orliczof \Phi p$, $p \in \maxexpat \gaussdensity$ are
 equal. In turn, equality of spaces implies the equivalence norms. It
 is possible to derive explicit bounds by choosing a $t$ such that the
 right-hand-side is smaller or equal to 1, see
 \cite{siri|trivellato:2020-SPL}.

\section{Gaussian Orlicz-Sobolev spaces}
In this final section, I plan to extend the exponential statistical
bundle's construction to allow for differentiable densities in the
Gaussian space. Namely, I suggest a model with sub-exponential random
variables of the Gaussian space whose weak derivatives are
sub-Gaussian random variables. The presentation focusses on functional
analytic properties, see \cite{adams|fournier:2003} and
\cite{brezis:2011FASSPDE} as a general reference.

\emph{As a model of the statistical manifold, let us define the
  function space
\begin{equation}
  \label{eq:sobolev+}
  W^1 \orliczpof {\coshtwo} {1,2} \gaussdensity =  \setof{f\in \orliczpof
    {\coshtwo} 1 \gaussdensity}{\partial_i f 
    \in \orliczpof {\coshtwo} 2 \gaussdensity, i=1,\dots,n} \ ,
\end{equation}
where the weak partial derivative $\partial_j f$ exists if
\begin{equation*}
  \scalarat \gamma {\partial_j f}{\phi} = \scalarat \gamma
  {f}{\delta_j \phi} \quad \text{for all} \quad \phi \in C^1_0(\reals^n) \ .
\end{equation*}}

The above definition of weak derivative coincides this the usual
definition of derivative in the sense of Schwartz distributions
because $\phi \leftrightarrow \phi \cdot \gamma$ is a bijection of
$C_0^\infty(\reals^n)$ and
\begin{equation*}
  \scalarat \gaussdensity f {\delta_i \phi} = - \int f(x) \pderivby
  {x_i} 
  \left(\phi(x)\gaussdensity(x)\right) \ dx \ .
\end{equation*}

Weak derivatives do not provide tools for non-linear
computations. Hence we want to recall the relation between translation
and weak derivative, the weak version of Calculus's fundamental
theorem. Let be given a locally integrable real mapping
$G \in L^1_\text{loc}(\reals)$, and assume there exists a locally
integrable function $G'$ which is the weak derivative of $G$, that is,
\begin{equation}\label{eq:dG-weak}
  \int G(x) \phi'(x) \ dx = - \int G'(x) \phi(x) \ dx \ , \quad \phi
  \in C_0^\infty(\reals) \ .
\end{equation}
Define the translation $\tau_h G$, $h\in\reals$, by $t_h G(x) =
G(x-h)$, $h\in\reals$. It follows immediatly $\tau_h G \in
L^1_\text{loc}(\reals)$ and 
\begin{align*}
  \int (\tau_{-h} G(x) - G(x)) \phi(x) \ dx
  &= \int G(x+h) \phi(x) \
    dx - \int G(x) \phi(x) \
    dx \\
  &= \int G(x) ( \phi(x-h) - \phi(x)) \ dx \\
  &= \int G(x) \left. \phi(x-sh) \right|_{s=0}^{s=1} \ dx \\
  &= - h \int G(x) \int_0^1 \phi'(x-sh)\ ds \ dx \\
  &= - h \int_0^1 \int G(x) \phi'(x-sh) \ dx \ ds\\
  &= h \int_0^1  \int G'(x) \phi'(x-sh) \ dx \ ds = \\
  &= h \int_0^1  \int G'(x+sh) \phi'(x) \ dx \ ds = \\
  &= \int \left(h \int_0^1 G'(x+sh) \ ds\right) \phi(x) \ dx 
\end{align*}
As $\phi$ is any function in $C^\infty_0(\reals)$, we have proved that
\begin{equation}\label{eq:delta-G}
  \tau_{-h} G - G = h \int_0^1 \tau_{-sh} G' \ ds = h G' + h \int_0^1 (\tau_{-sh} G' - G') \ ds
\end{equation}
in $L^1_\text{loc}(\reals)$. In particular, if $G'$ is bounded by a
constant $K$, then $G$ is almost surely $K$-Lipschitz, $\avalof{G(x-h)
  - G(x)} \leq K \avalof h$.

Conversely, if \cref{eq:delta-G} holds, then \cref{eq:dG-weak} holds,
see, for example, \cite[Lemma~8.1-2]{brezis:2011FASSPDE}.

The argument above extends to $n$-variate functions $f \in
W^{1,1}_\text{loc}(\reals^n)$, that is, $f, \partial_i f
\in L^1_\text{loc}(\reals^n)$, $i = 1,\dots,n$, by considering, for
each $h \in \reals^n$, the univariate function $t \mapsto \tau_{th} f$
defined by $\tau_{th} f(x) = f(x-th)$. We obtain 
\begin{equation}\label{eq:delta-f}
  \tau_{-th} f - f = t \nabla f \cdot h + t \int_0^1 (\tau_{-sth}
  \nabla f - \nabla f) \cdot h \ ds \ , 
\end{equation}
where the equality holds in $L^1_\text{loc}(\reals^n)$. The same
equality holds in all function space whose elements are locally
integrable.

The result about the functional differentiability of translations is
the following. \emph{For all
  $f \in W^1 \orliczpof {\cosh_2} {1,2} \gaussdensity$ the following
  first increment equation holds,
\begin{equation*}
  (\tau_{-th} f - f) - t \nabla f \cdot h = t \int_0^1 (\tau_{-sth}
  \nabla f - \nabla f) \cdot h \ ds \ ,
\end{equation*}
  and differentiability holds in $L^{\infty-0}(\gaussdensity) = \cap_{\alpha>1} L^\alpha(\gamma)$.
} In fact, the translations are continous in all $L^\alpha(\gaussdensity)$.

It is possible to exted the previous result to a property of the
derivative of the composite function $G\circ f$.  The increment of the
composition expands as
\begin{multline}\label{eq:delta-Gf}
  G(f(x + th)) - G(f(x)) = G(f(x) + (f(x+th) - f(x)) = \\
  (f(x+th) - f(x)) G'(f(x)) + \\ (f(x+th) - f(x)) \int_0^1 (G'(f(x) + s
  (f(x+th)-f(x))) - G'(f(x))) \ ds \ ,
\end{multline}
and the weak derivative of the composite function exists if $G'$ is
bounded and $f \in W^1 \orliczpof {\coshtwo} {1,2}
\gaussdensity$. However, differentiability holds in
$L^{\infty-0}(\gaussdensity)$.

An interesting example of application is the neuron of a neural
network \cite[Ch. 18]{efron|hastie:2016}. If
$f_1, \dots, f_k \in W^1 \orliczpof {\coshtwo} {1,2} \gaussdensity$,
$G$ is an activation function with linear growth, for example, the
linear rectifier
$G(x) = x^+$, and $a_i, w_{ij}, b_i$ are given constants, then
\begin{equation*}
  \sum_{i=1}^h a_i G\left(\sum_{j=1}^k w_{ij} f_j - b_i\right) \in W^1 \orliczpof {\coshtwo} {1,2}
\gaussdensity \ .
\end{equation*}

A crucial feature of the space defined in \cref{eq:sobolev+} is
the fact that \emph{each element of $W^1 \orliczpof {\coshtwo} {1,2}
  \gaussdensity$ has a continuous version}. The embedding
\begin{equation*}
  \orliczpof {\cosh_2} 1 \gaussdensity, \orliczpof {\cosh_2} 2
\gaussdensity \subset \cap_{\alpha \geq 1} L^\alpha(\gaussdensity)
\end{equation*}
allows to use the standard Sobolev inequalities to our
case. $W_{\text{loc}}^{1,\alpha}(\reals^n)$ denotes the space of
funcions whose restriction to each open ball
$B_\rho = \setof{x \in \reals^n}{\avalof x < \rho}$ is
$\alpha$-integrable, together with all weak partial
derivatives. $C^\lambda(\overline B_\rho)$ denotes $\lambda$-H\" older
functions on the closed ball.

\emph{\begin{enumerate}
    \item
  The following restriction and embedding hold true and are continuous,
  \begin{equation*}
    W^1 \orliczpof {\cosh_2} {1,2} \gaussdensity \to
    W^{1,\alpha}(\overline B_\rho) \subset C^\lambda(\overline
    B_\rho)\ , \quad \rho > 0 \ , \quad 0 < \lambda < 1 \ . 
  \end{equation*}
\item
  The following inclusions hold true and are continuous:
  \begin{equation*}W^1 \orliczpof {\cosh_2} {1,2} \gaussdensity \subset \cap_{\alpha
    \geq 1} W^{1,\alpha}_{\textnormal{loc}} \cap \orliczpof
  {\cosh_2} 1 \gaussdensity \subset
  C(\reals^n) \cap \orliczpof {\cosh_2} 1 \gaussdensity \ ,
\end{equation*}
where the space of continuous functions $C(\reals)$ is endowed with
the uniform convergence on compact sets.
\end{enumerate}
}

The embedding are easily verified.  If
$f \in \orliczpof {\cosh-1} 1 \gaussdensity$, then, for all
$k \in \naturals$, the inequalities $x^{2k}/(2k)! \leq \cosh_2(x)$ and
$(2\pi)^{-n/2} \euler^{-\rho^2/2} \leq \gaussdensity(x)$ for
$x \in B_\rho$ imply the inequality
 \begin{multline*}
   \frac{(2\pi)^{-n/2}\euler^{-\rho^2/2}}{(2k)!} \int_{B_\rho}
   \left(\frac{f(x)}{\normat {\orliczpof {\cosh_2} 1 \gaussdensity} f}\right)^{2k} \
   dx \leq \\ \int \cosh_2\left(\frac{f(x)}{\normat
       {\orliczpof {\cosh_2} 1 \gaussdensity} f}\right) \gamma(x) \ dx \leq 1 \ ,
 \end{multline*}
so that
\begin{equation*}
  \normat {L^{2k}(B_\rho)} f ^{2k} \leq (2\pi)^{n/2} (2k)! \euler^{\rho^2/2}
  \normat {\orliczpof {\cosh_2} 1 \gaussdensity} f \ .
\end{equation*}
A similar argument applies to the weak partial derivatives. Now we can
use the Sobolev embedding theorem, see
\cite[Th. 4.12]{adams|fournier:2003}.

Let us conclude by explicitly reviewing the main properties of our
space. 

\emph{\begin{enumerate}
\item The space $W^1 \orliczpof {\cosh_2} {1,2}
\gaussdensity$ contains the constants and all polinomial up to order
2.
\item Each element has a continuous version.
\item If $G : \reals \to \reals$ is the primitive of a bounded
  function, then $G \circ f \in W^1 \orliczpof {\cosh_2} {1,2} \gaussdensity$.
\item $\min(f,g), \max(f,g) \in W^1 \orliczpof {\cosh_2} {1,2} \gaussdensity$.
\end{enumerate}
}

Moreover, our space is a Banach space: \emph{The mapping
  \begin{equation*}
  W^1 \orliczpof \Phi {1,2} \gaussdensity \ni f \mapsto \normat
  {\orliczpof \Phi 1 \gaussdensity} f + \sum_{i=1}^n \normat {\orliczpof \Phi 2
    \gaussdensity} {\partial_i f}  
  \end{equation*}
is a complete norm and thus defines a Banach space.}
The argument is a standard one in Functional Analysis. The weak gradient $\nabla$ is a closed operator from $\orliczpof \Phi
1 \gaussdensity \to \left(\orliczpof \Phi 2 \gaussdensity\right)^n$,
that is, the graph of $\nabla$ is closed in $\orliczpof \Phi
1 \gaussdensity \times \left(\orliczpof \Phi 2
  \gaussdensity\right)^n$. In fact, given a converging sequence in the
graph, say $f_n \to f$ and $\partial_i f
\to f^i$, it holds
\begin{equation*}
  \scalarat \gaussdensity {\partial_i f} {\phi} = \scalarat \gaussdensity {f}
  {\delta_i \phi} = \lim_n \scalarat \gaussdensity {f_n} {\delta_i \phi}
  = \lim_n \scalarat \gaussdensity
  {\partial_i f} {\phi} = \scalarat \gaussdensity {f^i} {\phi} \ .
\end{equation*}
The identification of the space $W^1 \orliczpof \Phi {1,2}
\gaussdensity$ with the graph of $\nabla$ provides a complete
norm.

Finally, all exponential Orlicz spaces are isomorphic in a
maximal exponential model. This, in turn implies the isomorphism
\begin{equation*}
  W^1 \orliczpof {\coshtwo} {1,2}
\gaussdensity \leftrightarrow   W^1 \orliczpof {\coshtwo} {1,2}
{p \cdot \gaussdensity} \quad p \in \maxexpat \gaussdensity \ . 
\end{equation*}
It follows that this space is suitable as a model of the fibers of a
statistical bundle.

\section{Selected bibliography}
\label{sec:select-bibl}
General monograph on IG are \cite{kass|vos:1997},
\cite{amari|nagaoka:2000}, and
\cite{Ay|Jost|Le|Schwachhofer:2017IGbook}.  The type of analytic
framework I suggest to use in IG is that of Banach manifolds as in
\cite{lang:1995}. A first version of this project has been developed
in
\cite{pistone|sempi:95,gibilisco|pistone:98,pistone|rogantin:99,cena:2002}. In
a series of papers \cite{pistone:2013Entropy},
\cite{lods|pistone:2015}, \cite{pistone:2017-GSI2017}, and
\cite{pistone:2018-IGAIA-IV}, we have explored a version of the
non-parametric Information Geometry (IG) for smooth densities on
$\reals^n$. Especially, we have considered the IG associated to Orlicz
spaces on the Gaussian space. The analysis of the Gaussian space is
discussed, for example, in \cite{malliavin:1997}, and
\cite{nourdin|peccati:2012}.

This set-up provides a way to construct a statistical manifold
modelled on functional spaces of smooth densities, but other modelling
options are in fact available, for example the global analysis methods
of \cite{kriegl|michor:1997} and the approach based on deformed
exponentials \cite{newton:2019LNCS}. For the two examples involving
derivatives, see \cite{hyvarinen:2005}, \cite{lods|pistone:2015},
\cite{otto:2001}, and \cite{lott:2008calculations}.

 We have worked with a restricted type of Young functions. See the
 more general cases in \cite[Ch.~II]{musielak:1983} and
 \cite[Ch.~VII]{adams|fournier:2003}. There is a large literature
 about sub-exponential random variables, for example,
 \cite{buldygin|kozachenko:2000}, \cite{vershynin:2018-HDP},
 \cite{wainwright:2019-HDS}. Application to IG is discussed in
 \cite{siri|trivellato:2020-SPL}. The need to control the product of
 two random variables in $\orliczof {(\cosh-1)} \mu$ appears, notably,
 in the study of the covariant derivatives of the statistical bundle,
 see \cite{gibilisco|pistone:98}, \cite{lott:2008calculations},
 \cite{pistone:2018Lagrange}, and
 \cite{chirco|malago|pistone:2020-2009.09431}.

See \cite{pistone:2013GSI}, and \cite{pistone:2018-IGAIA-IV} for all
details about the contruction of the statistical bundle that are
missing here. In particular the isomorphism of fibers is discussed in
detail in \cite{santacroce|siri|trivellato:2016} and
\cite{santacroce|siri|trivellato:2018}. For Malliavin's calculus, see
\cite{malliavin:1997} and \cite{malliavin:1995}. General references
about the functional background used in the construction of the
Gaussian Orlicz-Sobolev space are \cite{adams|fournier:2003}, and
\cite{brezis:2011FASSPDE}.

\subsection*{acknowledgements}The author acknowledges the support by de Castro
  Statistics, Collegio Carlo Alberto, Turin, Italy. He is a member of
  GNAMPA-INDAM. The author thanks an anonymous reviewer and L. Malag\`o for their helpful comments.

\bibliographystyle{splncs04}

\end{document}